\documentclass[twoside]{article}
\usepackage{amsmath, amssymb,graphicx,wrapfig, url}
\title{The Combinatorics of  Flat
Folds: a Survey\footnote{This article appears in {\em Origami$^3$: The Third
International Meeting of Origami Science, Mathematics, and Education} published by AK Peters,
2002.}}

\author{Thomas C. Hull\footnote{Department of Mathematics, Western New England University, Springfield, MA 01119, thull@wne.edu}}

\date{}

\newtheorem{thm}{Theorem}[section]

\begin{document}

\maketitle

\begin{abstract}
We survey results on the foldability of flat origami models. The main topics are the question
of when a given crease pattern can fold flat, the combinatorics of mountain and valley
creases, and counting how many ways a given crease pattern can be folded. In particular, we
explore generalizations of Maekawa's and Kawasaki's Theorems, develop a necessary and
sufficient condition for a given assignment of mountains and valleys to fold up in a special
case of single vertex folds, and describe recursive formulas to enumerate the number of ways
that single vertex in a crease pattern can be folded.
\end{abstract}

\section{Introduction}

It is safe to say that in the study of the mathematics of origami, flat
origami has received the most attention.  To put it simply, a {\em flat origami model} is
one which can be pressed in a book without (in theory) introducing new creases.  We say
``in theory" because when one actually folds a flat origami model, slight errors in
folding will often make the model slightly non-flat.  In our analysis, however, we
ignore such errors and assume all of our models are perfectly folded. 
We also assume that our paper has zero thickness and that our creases have no width.  
It is surprising how rich the results are using a purely combinatorial
analysis of flat origami.  In this paper we introduce the basics of this approach,
survey the known results, and briefly describe where future work might lie. 

\begin{figure}
\centerline{\includegraphics[scale=.5]{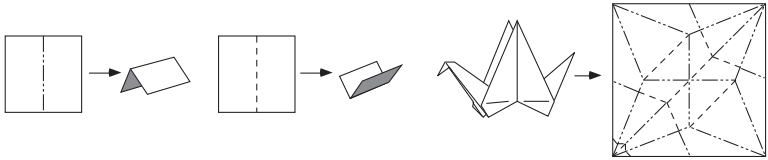}}
\caption{Mountain creases, valley creases, and the crease pattern for the flapping bird
with MV assignment shown.}\label{hufig1}
\end{figure}

First, some basic definitions are in order.  A {\em fold}
refers to any folded paper object, independent of the number of folds done
in sequence.  The {\em crease
pattern} of a fold is a planar embedding of a graph which represents the creases that
are used in the final folded
object.  (This can be thought of as a structural blueprint of the fold.)
Creases come in two types: {\em mountain
creases}, which are convex, and {\em valley creases}, which are concave (see Figure
\ref{hufig1}). Clearly the type of a crease depends on which side of the paper we look at,
and so we assume we are always looking at the same side of the paper. 

We also define a {\em mountain-valley (MV) assignment} to be a function mapping the set
of all creases to the set $\{M,V\}$.  In other words, we label each crease mountain or
valley.   MV assignments that can actually be folded are called {\em valid}, while those which
do not admit a flat folding (i.e. force the paper to self-intersect in some way) are called
{\em invalid}.

There are two basic questions on which flat-folding research has focused:

\newcounter{hull1}
\begin{list}
{\arabic{hull1}.}{\usecounter{hull1}
\setlength{\parsep}{0in}\setlength{\itemsep}{0in}}
\item Given a crease pattern, without an MV assignment, can we tell whether it can flat fold?
\item If an MV assignment is given as well, can we tell whether it is valid?
\end{list}

These are also the focus of this survey. We will not discuss the special cases of flat origami
tesselations, origami model design, or other science applications.

\section{Classic single vertex results}

We start with the simplest case for flat origami folds. We define a {\em single vertex fold}
to be a crease pattern (no MV assignment) with only one vertex in the interior of the paper and
all crease lines incident to it. Intersections of creases on the boundary of the paper clearly
follow different rules, and nothing of interest has been found to say about them thus far
(except in origami design; see \cite{hulan1}, \cite{hulan2}).
 A single vertex fold which is
known to fold flat is called a {\em flat vertex fold}. We present a few basic theorems
relating to necessary and sufficient conditions for flat-foldability of single vertex
folds. These theorems appear in their cited references without proof.  While Kawasaki,
Maekawa, and Justin undoubtedly had proofs of their own, the proofs presented below appear in
\cite{huhul1}.

\begin{thm}[Kawasaki \cite{hukaw2}, Justin \cite{hujus1},
\cite{hujus2}]\label{hukj}
Let $v$ be a vertex of degree $2n$ in a single vertex fold
and let
$\alpha_1, ..., \alpha_{2n}$ be the consecutive angles between the creases.  Then  $v$ is
a flat vertex fold if and only if
\begin{equation}\label{huiso}
\alpha_1-\alpha_2+\alpha_3-\cdots -\alpha_{2n}=0.
\end{equation}
\end{thm}

\noindent{\bf Proof:} Consider a simple closed curve
which winds around the vertex.  This curve mimics the path of an ant
walking around the vertex on the surface of the paper after it is folded.  We
measure the
angles the ant crosses as positive when traveling to the left and negative when walking to
the right.  Arriving at the point where the ant started means that this alternating sum is
zero. The converse is left to the reader; see
\cite{huhul1} for more details. $\Box$

\begin{thm}[Maekawa, Justin \cite{hujus2}, \cite{hukas}]\label{humj}
Let $M$ be the number of mountain creases and $V$ be the number of valley
creases adjacent to a vertex in a single vertex fold.  Then $M-V=\pm 2$.
\end{thm}

\noindent{\bf Proof:} (Siwanowicz) If $n$ is the number of creases, then $n=M+V$. 
Fold the paper flat and 
consider the cross-section
obtained by clipping the area near the vertex from the paper; the cross-section
forms a flat polygon.  If we view each interior
$0^\circ$ angle as a valley crease and each interior $360^\circ$ angle as a mountain
crease, then $0V+360M=(n-2)180=(M+V-2)180$, which gives $M-V=-2$. On the other hand,
if we view each
$0^\circ$ angle as a \emph{mountain} crease and each $360^\circ$ angle as a
\emph{valley} crease (this corresponds to flipping the paper over), then
we get $M-V=2$.  $\Box$

\vspace{.1in}

In the literature, Theorem \ref{hukj} and \ref{humj} are referred to as Kawasaki's Theorem
and Maekawa's Theorem, respectively.  
Justin \cite{hujus3}
refers to equation (\ref{huiso}) as the {\em isometries condition}. 
Kawasaki's Theorem is sometimes
stated in the equivalent form that the sum of alternate angles around $v$ equals
$180^\circ$, but this is only true if the vertex is on a flat sheet of paper. 
Indeed, notice that the proofs of the Kawasaki's and Maekawa's Theorems do not use the fact
that
$\sum \alpha_i = 360^\circ$.  Thus these two theorems are also valid for single vertex
folds where $v$ is at the apex of a cone-shaped piece of paper.  We will require
this generalization in sections 4 and 5.

Note that while Kawasaki's Theorem assumes that the vertex has even degree, Maekawa's
Theorem does not.  Indeed, Maekawa's Theorem can be used to prove this fact.  Let
$v$ be a single vertex fold that folds flat and let $n$ be the degree of $v$.
Then $n=M+V=M-V+2V=\pm 2 + 2V$, which is even.

\section{Generalizing Kawasaki's Theorem}

\begin{figure}
\centerline{\includegraphics[scale=.5]{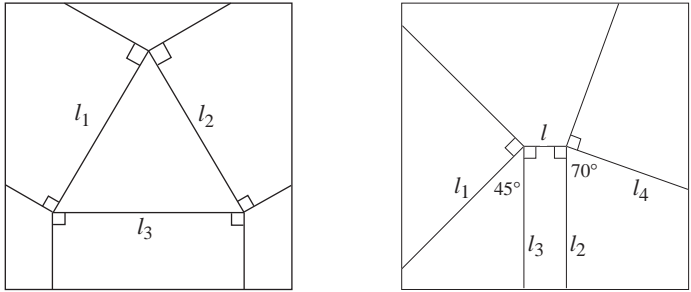}}
\caption{Two impossible-to-fold-flat folds.}\label{hufig2}
\end{figure}

Kawasaki's Theorem gives us a complete description of when a single vertex in a crease
pattern will (locally) fold flat.   Figure \ref{hufig2} shows two examples of crease patterns
which satisfy Kawasaki's Theorem at each vertex, but which will not fold flat.  The example on
the left is from \cite{huhul1}, and a simple argument shows that no two of the creases $l_1,
l_2, l_3$ can have the same MV parity.  Thus no valid MV assignment for
the lines
$l_1, l_2, l_3$ is possible.  The example on the right has valid MV assignments, but still
fails to fold flat.  The reader is encouraged to copy this crease
pattern and try to fold it flat, which will reveal that some flap of paper will have to
intersect one of the creases.  However, if the location of the two
vertices is changed relative to the border of the paper, or if the crease $l$ is made longer,
then the crease pattern {\em will} fold flat.  

This illustrates how difficult the question of flat-foldability is for
multiple vertex folds.  Indeed, in 1996 Bern and Hayes \cite{huber} proved
that the general question of whether or not a given crease pattern can fold flat is
NP-hard.  Thus one would not expect to find easy necessary and sufficient conditions for
general flat-foldability.  

We will present two efforts to describe general flat-foldability.  The first has to do
with the realization that when we fold flat along a crease, one part of the paper is being
reflected along the crease line to the other side.  Let us denote $R(l_i)$ to be the
reflection in the plane, $\mathbb R^2$, along a line $l_i$.  

\begin{thm}[Kawasaki \cite{hukaw0}, \cite{hukaw3}, Justin \cite{hujus1},
\cite{hujus3}]\label{hureflect}
Given a multiple vertex fold, let $\gamma$ be any closed, vertex-avoiding curve
drawn on the crease pattern which crosses crease lines $l_1,...$, $l_{n}$, in order.  Then,
if the crease pattern can fold flat, we will have 
\begin{equation}\label{huref}
R(l_1)R(l_2)\cdots R(l_n)=I
\end{equation}
where $I$ denotes the identity transformation.
\end{thm}

Although a rigorous proof of Theorem \ref{hureflect} does not appear in the literature we
sketch here a proof by induction on the number of vertices.  In the base case, we are given a
single vertex fold, and it is a fun exercise to show that condition (\ref{huref}) is
equivalent to equation (\ref{huiso}) in Kawasaki's Theorem (use the fact that the composition
of two reflections is a rotation).  The induction step then proceeds by breaking the curve
$\gamma$ containing $k$ vertices into two closed curves, one containing $k-1$ vertices and one
containing a single vertex (the
$k$th).  

The condition (\ref{huref}) is not a sufficient condition for
flat-foldability (the crease patterns in Figure \ref{hufig2} are counterexamples here as
well).  In fact, as the induction proof illustrates, Theorem \ref{hureflect} 
extends Kawasaki's Theorem to as general a result as possible.  

In \cite{hujus3} Justin proposes a necessary and sufficient condition for general
flat-foldability, although as Bern and Hayes predicted, it is not very computationally
feasible. To summarize, let $C$ be a crease pattern for a flat origami model, but for the
moment we are considering the boundary of the paper as part of the graph.  If $E$ denotes the
set of edges in $C$ embedded in the plane, then we call $\mu(E)$ the {\em f-net}, which is the
image of all creases and boundary of the paper after the model has been folded.  We then call
$\mu^{-1}(\mu(E))$ the {\em s-net}.  This is equivalent to
imagining that we fold carbon-sensitive paper, rub all the
crease lines firmly, and then unfold.  The result will be the $s$-net.

Justin's idea is as follows: Take all the faces of the $s$-net which get mapped by $\mu$ to the
same region of the $f$-net and assign a {\em superposition order} to them in accordance to
their layering in the final folded model.  One can thus try to fold a given crease pattern by
cutting the piece of paper along the creases of the $s$-net, tranforming them under
$\mu$, applying the superposition order, and then attempting to glue the paper back together. 
Justin describes a set of three intuitive {\em crossing conditions} (see \cite{hujus3}) which
must not happen along the $s$-net creases during the glueing process if the model is to be
flat-foldable -- if this can be done, we say that the {\em non-crossing condition} is
satisfied.  Essentially Justin conjectures that a crease pattern folds flat if and only if the
non-crossing condition holds. 
Although the spirit of this approach seems to accurately reflect the flat-foldability of
multiple vertex folds, no rigorous proof appears in the literature; it seem that this is an
open problem.

\section{Generalizing Maekawa's Theorem}

To extend Maekawa's Theorem to more than one vertex, we define interior vertices in a flat
multiple vertex fold to be {\em up vertices} and {\em down vertices} if they locally have
$M-V=2$ or
$-2$, respectively.  We define a crease line to be an {\em interior} crease if its
endpoints lie in the interior of the paper (as opposed to on the boundary), and
consider any crease line with both endpoints on the boundary of the paper to
actually be two crease lines with an interior vertex of degree 2 separating them.

\begin{thm}[Hull \cite{huhul1}]
Given a multiple vertex flat fold, let $M$ (resp. $V$) denote the number of mountain (resp. 
valley) creases, $U$ (resp. $D$) denote the number of up (resp. down) vertices, and $M_i$
(resp. $V_i$) denote the number of interior mountain  (resp. valley) creases.  Then
$$M-V=2U-2D-M_i+V_i.$$
\end{thm}

Another interesting way to generalize Maekawa's Theorem is to explore restrictions which turn
it into a sufficiency condition.  In the case where all of the angles
around a single vertex are equal, an MV assignment with $M-V=\pm 2$ is guaranteed to be valid. 
This observation can be generalized to sequences of consecutive equal angles around a
vertex.

Let us denote a single vertex fold by
$v=(\alpha_1,... \alpha_{2n})$ where the $\alpha_i$ are consecutive angles between the
crease lines.  We let
$l_1,..., l_{2n}$ denote the creases adjacent to 
$v$ where $\alpha_i$
is the angle between creases $l_i$ and $l_{i+1}$ ($\alpha_{2n}$
is between $l_{2n}$ and $l_1$).  

 If $l_i,...,l_{i+k}$
are consecutive crease lines in a single vertex fold which have been given a MV
assignment, let $M_{i,...,i+k}$ be the number of mountains and $V_{i,...,i+k}$ bethe
number of valleys among these crease lines. We say that a given MV assignment is valid for the
crease lines
$l_i, ..., l_{i+k}$  if the MV assignment can be
followed to fold up these crease lines without forcing the paper to self-intersect.  (Unless
these lines include all the creases at the vertex, the result will be a cone.) The necessity
porttion oof the following result appears in \cite{huhul2}, while sufficiency is new.

\begin{thm}\label{hullth1}
Let $v=(\alpha_1,...,\alpha_{2n})$ be a single vertex fold in either a piece of paper or a
cone, and suppose we have
$\alpha_i= \alpha_{i+1}=\alpha_{i+2}=\cdots =\alpha_{i+k}$  for some
$i$ and $k$.  Then a given MV assignment is valid for $l_i,..., l_{i+k+1}$ if
and only if 
$$M_{i,...,i+k+1}-V_{i,...,i+k+1} = \left\{\begin{array}{cl}
0 & \mbox{when $k$ is even}\\
\pm 1 & \mbox{when $k$ is odd.}\end{array}\right.$$
\end{thm}

\begin{figure}[h]
\centerline{\includegraphics[scale=.5]{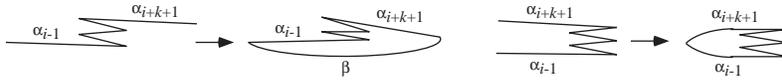}}
\caption{Applying Maekawa when $k$ is even (left) and odd (right).}
\label{hufig3}
\end{figure}

\noindent{\bf Proof:} Necessity follows by applications of Maekawa's
Theorem.  If $k$ is even, then the cross-section of the paper around the
creases in question might look as shown in the left of Figure \ref{hufig3}. If we consider only
this sequence of angles and imagine adding a section of paper with angle $\beta$ to
connect the loose ends at the left and right (see Figure \ref{hufig3}, left), then we'll have
a flat-folded cone which must satisfy Maekawa's
 Theorem.  The angle $\beta$ adds two extra creases, both of which
must be mountains (or valleys).  We may assume that the vertex points up, and thus
we subtract two from the result of Maekawa's Theorem to get
$M_{i,...,i+k+1}-V_{i,...,i+k+1}=0$.

If $k$ is odd (Figure \ref{hufig3}, right), then this angle sequence, if considered by
itself, will have the loose ends from angles $\alpha_{i-1}$ and $\alpha_{i+k+1}$ pointing
in the same direction.  If we glue these together (extending them if necessary) 
then Maekawa's Theorem may be applied.  After subtracting (or adding) one to the result of
Maekawa's Theorem because of the extra crease made when gluing the loose flaps, we get
$M_{i,...,i+k+1}-V_{i,...,i+k+1}=\pm 1$. 

For sufficiency, we proceed by induction on $k$.  The result is trivial for the base cases
$k=0$ (only one angle, and the two neighboring creases will either be M, V or V, M) and
$k=1$ (two angles, and all three possible ways to assign 2 M's and 1 V, or vice-versa, can
be readily checked to be foldable).  For arbitrary $k$, we will always be able to find
two adjacent creases $l_{i+j}$ and $l_{i+j+1}$ to which the MV assignment assigns
opposite parity. Let $l_{i+j}$ be M and $l_{i+j+1}$ be V.  We make these folds and we
can imagine that $\alpha_{i+j-1}$ and $\alpha_{i+j}$ have been fused into the other layers of
paper, i.e. removed.  The value of
$M-V$ will not have changed for the remaining sequence $l_i, ..., l_{i+j-1}, l_{i+j+2}, ...,
l_{i+k}$ of creases, which are flat-foldable by the induction hypothesis. $\Box$

\section{Counting valid MV assignments}

We now turn to the question of counting how many different ways we can fold a flat origami
model.  By this we mean, given a crease pattern that is known to fold flat, how many
different valid MV assignments are possible?  

We start with the single vertex case.  Let $C(\alpha_1,...,\alpha_{2n})$ denote the
number of valid MV assignments possible for the vertex fold $v=(\alpha_1,...,\alpha_{2n})$.

\begin{figure}[h]
\centerline{\includegraphics[scale=.6]{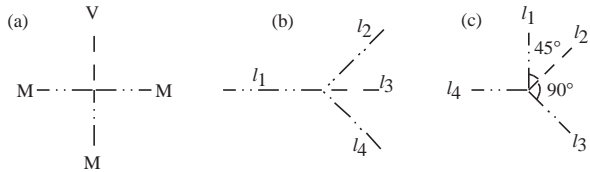}}
\caption{The three scenarios for vertices of degree 4.}\label{hufig4}
\end{figure}

An an example, consider the case where $n=2$ (so we have 4 crease lines at $v$).  We compute
$C(\alpha_1,\alpha_2,\alpha_3,\alpha_4)$ using Maekawa's Theorem. Its value will depend on the
type of symmetry the vertex has, and the three possible situations are depicted in Figure
\ref{hufig4}.  $C(90,90,90,90)=8$ because any crease could be the
``odd creasee out" and the vertex could be up or down.  In Figure
\ref{hufig4} (b) we have only mirror symmetry, and by Theorem \ref{hullth1},
$M_{2,3,4}-V_{2,3,4}=\pm 1$.  Thus
$l_2,l_3,l_4$ must have 2 M's and 1 V or vice versa; this determines $l_1$'s parity, giving
$C(\alpha_1,...,\alpha_4)=6$. In Figure \ref{hufig4} (c)
$M_{1,2}-V_{1,2}=0$, so $l_1$ and $l_2$ can be M,V or V,M, and the other two must be
both M or both V, giving $C(\alpha_1,...,\alpha_4)=4$.

The example in Figure \ref{hufig4}(a) represents the case with no restrictions.
This appears whenever all the angles are equal around $v$, giving
$C(\alpha_1,...,\alpha_{2n})$ $= 2{2n\choose n-1}$.  The idea in Figure
\ref{hufig4} (c), where we pick the smallest angle we see and let its creases be M,V or
V,M, can be applied inductively to give the lower bound in the following (see
\cite{huhul2} for a full proof):

\begin{thm}\label{hullth2}
Let $v=(\alpha_1,...,\alpha_{2n})$ be the vertex in a flat vertex fold, on either a
flat piece of paper or a cone.  Then
$$2^n\leq C(\alpha_1,...,\alpha_{2n})\leq 2{2n\choose n-1}$$
are sharp bounds.
\end{thm}

A formula for $C(\alpha_1,...,\alpha_{2n})$ seems out of reach, but using the
equal-angles-in-a-row concept, recursive formulas exist to compute this quantity in
linear time.

\begin{thm}[Hull, \cite{huhul2}]\label{hullth3}
Let $v=(\alpha_1,...,\alpha_{2n})$ be a flat vertex fold in either a piece of paper 
or a cone, and suppose we have
$\alpha_i= \alpha_{i+1}=\alpha_{i+2}=\cdots =\alpha_{i+k}$ and $\alpha_{i-1}>
\alpha_i$ and
$\alpha_{i+k+1}>\alpha_{i+k}$ for some
$i$ and $k$.  Then
$$C(\alpha_1,...,\alpha_{2n}) = 
{k+2\choose
\frac{k+2}{2}}C(\alpha_1,...,\alpha_{i-2},\alpha_{i-1}-\alpha_i+\alpha_{i+k+1},
\alpha_{i+k+2},..., \alpha_{2n})$$
if $k$ is even, and
$$C(\alpha_1,...,\alpha_{2n}) =
{k+2\choose \frac{k+1}{2}}C(\alpha_1,...,\alpha_{i-1},\alpha_{i+k+1}, ...,
\alpha_{2n})$$
if $k$ is odd.
\end{thm}

Theorem \ref{hullth3} was first stated in \cite{huhul2}, but the basic ideas behind it are
discussed by Justin in \cite{hujus3}.  

As an example, consider $C(20,10,40,50,60,60,60,60)$. 
Theorem \ref{hullth2} tells us that this qualtity lies between 16 and 112.  But using Theorem
\ref{hullth3} we see that $C(20,10,40,50,60,60,60,60)=$ ${2\choose 1}C(50, 50, 60,60,60,60)$ 
$= {2\choose 1}{3\choose 1}C(60,60,60,60)$
$= {2\choose 1}{3\choose 1}2{4\choose 1} =48$.

Not much is known about counting valid MV assignments for flat multiple vertex folds. 
While there are similarities with work done on counting the number of ways to fold up a grid of
postage stamps (see \cite{hukoe}, \cite{hulun1}, \cite{hulun2}), the questions asked are
slightly different.  For other work, see \cite{hujus3} and
\cite{huhul2}.

\section{Conclusion}

In conclusion, the results for flat-foldability seem to have almost completely exhausted the
single vertex case.  Open problems exist, however, in terms of global flat-foldability, and
very little is known about enumerating valid MV assignments for multiple vertex
crease patterns.


\begin{thebibliography}{99}


\bibitem{huber} Bern, M. and Hayes, B.,
``The complexity of flat origami",
{\em Proceedings of the 7th Annual ACM-SIAM Symposium on Discrete
  Algorithms}, (1996) 175--183.

\bibitem{hueh} Ewins, B. and Hull, T., personal communication, 1994.

\bibitem{huhul1} Hull, T., ``On the mathematics of flat origamis", {\em Congressus
Numerantium}, 100 (1994) 215-224.

\bibitem{huhul2} Hull, T., ``Counting mountain-valley assignments for flat folds",
{\em Ars Combinatoria}, to appear.

\bibitem{hujus1} Justin, J., ``Aspects mathematiques du pliage de papier" (in French),
in: H. Huzita ed., {\em Proceedings of
the First International Meeting of Origami Science and Technology},
Ferrara, (1989) 263-277.

\bibitem{hujus2} Justin, J., ``Mathematics of origami, part 9", {\em British Origami}
(June 1986) 28-30.

\bibitem{hujus3} Justin, J., ``Toward a mathematical theory of origami", in: K.
Miura ed., {\em Origami Science and Art: Proceedings of the Second International
Meeting of Origami Science and Scientific Origami}, Seian University of Art and
Design, Otsu, (1997) 15-29.

\bibitem{hukas} Kasahara, K. and Takahama, T., {\em Origami for the Connoisseur},
Japan Publications, New York, (1987).

\bibitem{hukaw0} Kawasaki, T. and Yoshida, M., ``Crystallographic flat origamis",
{\em Memoirs of the Faculty of Science, Kyushu University, Series A},
Vol. 42, No. 2 (1988), 153-157. 

\bibitem{hukaw2} Kawasaki, T., ``On the relation between mountain-creases and
valley-creases of a flat origami" (abridged English translation), in: H. Huzita ed., 
{\em Proceedings of
the First International Meeting of Origami Science and Technology},
Ferrara, (1989) 229-237.

\bibitem{hukaw1} Kawasaki, T., ``On the relation between mountain-creases and
valley-creases of a flat origami" (unabridged, in Japanese), Sasebo College of
Technology Reports, 27, (1990) 55-80.


\bibitem{hukaw3} Kawasaki, T., ``$R(\gamma)=I$", in: K.
Miura ed., {\em Origami Science and Art: Proceedings of the Second International Meeting
of Origami Science and Scientific Origami}, Seian University of Art and Design, Otsu,
(1997) 31-40.

\bibitem{hukoe} Koehler, J., ``Folding a strip of stamps", {\em Journal of
Combinatorial Theory}, 5 (1968) 135-152.

\bibitem{hulan1} Lang, R.J.,
``A computational algorithm for origami design", {\em Proceedings of the 12th Annual ACM
Symposium on Computational
  Geometry}, (1996), 98-105.

\bibitem{hulan2} Lang, R.J.,
{\em {TreeMaker} 4.0: A Program for Origami Design}, (1998)
\url{http://origami.kvi.nl/programs/TreeMaker/trmkr40.pdf}

\bibitem{hulun1} Lunnon, W.F., ``A map-folding problem", {\em Mathematics of
Computation}, 22, No. 101 (1968) 193-199.

\bibitem{hulun2} Lunnon, W.F., ``Multi-dimensional map folding", {\em The Computer
Journal}, 14, No. 1 (1971) 75-80.

\end{thebibliography}
\end{document}